\documentclass[a4paper,reqno,12pt]{amsart}

\usepackage{amsmath,amsthm,amssymb,amsfonts}
\usepackage{amscd,tikz-cd}
\usepackage[all]{xy} 
\usepackage{comment,todonotes}

\usepackage{mathrsfs}

\theoremstyle{plain} 
\newtheorem{thm}{Theorem}[section]
\newtheorem*{thm*}{Theorem}
\newtheorem{lem}[thm]{Lemma}
\newtheorem{prop}[thm]{Proposition}
\newtheorem{cor}[thm]{Corollary}

\newtheorem{athm}{Theorem}[section]

\newtheorem{acor}[athm]{Corollary}

\theoremstyle{definition}
\newtheorem{dfn}[thm]{Definition}
\newtheorem{eg}[thm]{Example}

\newtheorem{rem}[thm]{Remark}
\newtheorem{conj}[thm]{Conjecture}
\newtheorem*{conj*}{Conjecture}

\newtheorem*{ack}{Acknowledgements}

\theoremstyle{remark}
\newtheorem*{pf}{Proof}

\numberwithin{equation}{section}

\newcommand{\ZZ}{\mathbb{Z}}
\newcommand{\QQ}{\mathbb{Q}}
\newcommand{\RR}{\mathbb{R}}
\newcommand{\CC}{\mathbb{C}}
\newcommand{\PP}{\mathbb{P}}
\newcommand{\HH}{\mathbb{H}}

\newcommand{\BB}{\mathbb{B}}

\newcommand{\EE}{\mathbb{E}}

\newcommand{\UU}{\mathbb{U}}

\newcommand{\B}{\mathcal{B}}
\newcommand{\C}{\mathcal{C}}

\newcommand{\N}{\mathcal{N}}

\def\S{\mathcal{S}}

\DeclareMathOperator{\Aut}{Aut}

\DeclareMathOperator{\NS}{NS}
\DeclareMathOperator{\Trans}{T}
\DeclareMathOperator{\Pic}{Pic}
\DeclareMathOperator{\Amp}{Amp}
\DeclareMathOperator{\Nef}{Nef}
\DeclareMathOperator{\Eff}{Eff}
\DeclareMathOperator{\Mov}{Mov}
\DeclareMathOperator{\Hilb}{Hilb}

\DeclareMathOperator{\vcd}{vcd}
\DeclareMathOperator{\MW}{MW}

\DeclareMathOperator{\rank}{rank}

\DeclareMathOperator{\Isom}{Isom}
\DeclareMathOperator{\Bir}{Bir}

\begin{document}
\title
{Geometrical finiteness for automorphism groups via cone conjecture}

\author{Kohei Kikuta}

\address{Department of Mathematics, Graduate School of Science, Osaka University, Toyonaka Osaka, 560-0043, Japan.
\vspace{1mm}
\newline
\hspace*{4mm}
School of Mathematics, 
The University of Edinburgh, 
James Clerk Maxwell Building, 
Peter Guthrie Tait Road, Edinburgh, Scotland, EH9 3FD, UK.}
\email{kikuta@math.sci.osaka-u.ac.jp}

\begin{abstract}
This paper aims to establish
the geometrical finiteness for the natural isometric actions of (birational) automorphism groups on the hyperbolic spaces
for K3 surfaces, Enriques surfaces, Coble surfaces, and irreducible symplectic varieties. 

As an application, it follows that 
such groups are non-positively curved: relatively hyperbolic and ${\rm CAT(0)}$. 
In the case of K3 surfaces, 
we clarify the relationship between Kleinian lattices and $(-2)$-curves, and between convex-cocompact Kleinian groups and genus-one fibrations. 

\end{abstract}

\maketitle
\date{}
\markboth{KOHEI KIKUTA}{GEOMETRICAL FINITENESS}


\section{Introduction}
\subsection{Geometrical finiteness}

Geometrical finiteness is one of the central notions in the study of Kleinian groups 
or more generally, 
discrete subgroups of semisimple Lie groups of real rank one. 
This term describes associated quotient orbifolds such that 
all the interesting geometry goes on in some compact subset. 
Many people have contributed to the development of this idea, 
and among them, 
Bowditch provided 
several equivalent definitions of geometrical finiteness \cite{Bow1}. 
For further details, we also refer to \cite{Kap} and the references therein. 

\subsection{Main results}
The study of automorphism groups of algebraic varieties is a classical subject in algebraic geometry. 
For certain algebraic varieties $X$, 
the torsion-free part of the N\'{e}ron--Severi group admits a hyperbolic lattice structure, 
and the representation of the automorphism group has a finite kernel, 
which establishes a connection to hyperbolic geometry.
%
%
We denote the associated isometric action on the hyperbolic space by
\[
\Aut(X)\to\Isom(\HH^{\rho_X-1}), 
\]
where $\rho_X$ is the Picard rank of $X$. 

The main result of this paper is 
the geometrical finiteness for automorphism groups of such varieties
further satisfying the Morrison--Kawamata cone conjecture: 
the existence of a rational polyhedral fundamental domain for the action 
on the effective nef cone. 
%
%
We also have a similar result for birational automorphism groups. 

\newpage

\begin{athm}\label{thm-geom-fin}
Let $G_X$ be one of the following groups: 
\begin{enumerate}
\item[(1)]
the automorphism group $\Aut(X)$ of a K3 surface $X$ over a field of characteristic different from $2$. 

\item[(2)]
the automorphism group $\Aut(X)$ of a non-supersingular K3 surface $X$ over a field of characteristic $2$. 

\item[(3)]
the automorphism group $\Aut(X)$ of an Enriques surface $X$ over an algebraically closed field of characteristic different from $2$. 

\item[(4)]
the automorphism group $\Aut(X)$ of 
a Coble surface 
$X$ over an algebraically closed field of characteristic $0$. 

\item[(5)]
the automorphism group $\Aut(X)$ of a smooth irreducible symplectic variety $X$ over a field of characteristic $0$. 

\item[(6)]
the birational automorphism group $\Bir(X)$ of a smooth irreducible symplectic variety $X$ over a field of characteristic $0$. 


\end{enumerate}
Then the representation 
$G_X\to\Isom(\HH^{\rho_X-1})$
is geometrically finite. 
\end{athm}








Geometrical finiteness has historically played a central role in the development of geometric group theory, 
where infinite groups are studied via their isometric actions on metric spaces of non-positive curvature, such as Gromov hyperbolic spaces and {\rm CAT(0)} spaces. 
As an immediate consequence of classical results on geometrically finite Kleinian groups, we obtain the following: 

\begin{acor}[Corollary \ref{cor-alternative} and \ref{cor-cat(0)}]
\label{intro-thm-alternative}
Let $G_X$ be a group as in Theorem \ref{thm-geom-fin}. 
\begin{enumerate}
\item[(1)]
$G_X$ is either virtually abelian or non-elementary relatively hyperbolic. 
\item[(2)]
$G_X$ is {\rm CAT(0)}. 
\end{enumerate}
\end{acor}


The study of automorphism groups of K3 surfaces has been a central topic since Nikulin's seminal work \cite{Nikulin}, and it remains an active area of research.
Among such groups, virtually abelian automorphism groups, including finite groups as a special case, form a well-studied and relatively tractable class.

In contrast, geometrical finiteness provides a natural framework for studying broader classes of automorphism groups, extending beyond the virtually abelian setting.
As further applications of Theorem~\ref{thm-geom-fin}, 
we clarify, in the case of K3 surfaces, the respective relationships between (Kleinian) lattices and $(-2)$-curves, and between convex-cocompact Kleinian groups and genus-one fibrations: 

\begin{acor}[Theorem \ref{chara-lattice} and \ref{chara-conv-cocpt}]
\label{intro-thm-lattice-cocpt}
Let $X$ be a K3 surface as in Theorem \ref{thm-geom-fin}. 
\begin{enumerate}
\item[(1)]
$\Aut(X)$ is a lattice 
if and only if 
$X$ has no $(-2)$-curves. 
\item[(2)]
$\Aut(X)$ is convex-cocompact 
if and only if 
either
\begin{enumerate}
\item
$X$ admits no genus-one fibrations, or
\item 
$X$ admits a genus-one fibration, and the Mordell--Weil group 
of any genus-one fibration is finite. 
\end{enumerate}
\end{enumerate}
\end{acor}

We also obtain a dynamical characterization of relative hyperbolicity using entropy (Corollary \ref{cor-dyn}). 

\subsection{Related works}
\begin{enumerate}
\item[(1)]
Takatsu has independently proved the geometrical finiteness for automorphism groups of K3 surfaces over the complex number field $\CC$ in a different way \cite{Takatsu}.   
As an application, 
Takatsu proves that the virtual cohomological dimension of automorphism groups is determined by the covering dimension of the blown-up boundaries of their ample cones.  
%
He also provides an affirmative example of Mukai's conjecture using sphere packings.

\item[(2)]


The relationship between algebraic geometry and hyperbolic geometry was systematically studied by Totaro in \cite{Tot, Totaro-hyp-geom}. 
The idea of combining the cone conjecture with hyperbolic geometry already appears in \cite{Tot}. 
Totaro also provided a sketch of a proof of the {\rm CAT(0)} property for automorphism groups of K3 surfaces in \cite{Totaro-hyp-geom}. 
Kurnosov--Yasinsky later provided a more detailed sketch of a proof of the {\rm CAT(0)} property for automorphism groups and birational automorphism groups of smooth irreducible symplectic varieties over $\CC$ in \cite{KY}. 
Their argument is along the lines of the discussion in \cite[Section 12.4]{Rat}. 
A similar line of argument is also used in Section \ref{sec-gen-statement} of the present paper.



\end{enumerate}

\subsection{Organization of the paper}

Section~\ref{Preliminaries} recalls the basic facts on geometrically finite Kleinian groups from hyperbolic geometry and several cones arising in algebraic geometry.
In Section~\ref{pf-main-thm}, 
we prove the main result, Theorem~\ref{thm-geom-fin}. 
To this end, we establish a general statement in hyperbolic geometry, Theorem~\ref{thm-general-statement}, which forms the technical core of the paper.
The final section is devoted to applications of Theorem \ref{thm-geom-fin}, including Theorems \ref{intro-thm-alternative} and \ref{intro-thm-lattice-cocpt}.


\begin{ack}
The author would like to thank Koji Fujiwara, Tomohiro Fukaya, Keiji Oguiso, Shin-ichi Oguni, and Taiki Takatsu for valuable discussions and helpful comments.
The author is also indebted to Arend Bayer, the University of Edinburgh, and the Hausdorff Research Institute for Mathematics funded by 
the Deutsche Forschungsgemeinschaft (DFG, German Research Foundation) 
under Germany's Excellence Strategy – EXC-2047/1 – 390685813,
for their kind hospitality. 
This work is supported by JSPS Overseas Research Fellow and JSPS KAKENHI Grant Number 21K13780. 

\end{ack}



\section{Preliminaries}\label{Preliminaries}

\subsection{Geometrically finite representations}
We briefly recall the basics of hyperbolic geometry and Kleinian groups, 
and provide the definition of the geometrical finiteness. 
For further details, readers are referred to \cite{Rat}. 


A {\it lattice} is a finitely generated free abelian group endowed with an integral non-degenerate symmetric bilinear form. 
Let $L$ be a lattice of signature $(1,n)$. 
Fix one of the two connected components of $\left\{v\in L_\RR \mid v^2>0\right\}$
called a {\it positive cone} and denoted by $\C$. 
The {\it hyperboloid model} $\HH^n$ of the {\it hyperbolic space} is defined by 
\[
\HH^n:=\left\{x\in \C \mid (x,x)=1\right\}
\]
endowed with the metric $d$ determined by $\cosh d(x,y):=(x,y)$. 
%
The {\it boundary} of $\HH^n$ is defined by
$\partial\HH^n:=\left(\partial\C\backslash\{0\}\right)/\RR_{>0}$. 
The set of rational points on $\HH^n$ is given by 
\[
\HH^n(\QQ):=\left\{x\in \HH^n \mid \RR x=\RR v \text{ for some }v\in\C\cap L_\QQ \right\}.
\]
Similarly, the set $\partial\HH^n(\QQ)$ of rational points on $\partial\HH^n$ is the image of $\left(\partial \C\cap L_\QQ\right)\backslash\{0\}$
with respect to the quotient $\partial\C\backslash\{0\}\to\partial\HH^n$. 
We set
$\overline{\HH^n}:=\HH^n\cup \partial\HH^n$ and
$\overline{\HH^n}(\QQ):=\HH^n(\QQ)\cup \partial\HH^n(\QQ)$. 


In this paper, we also consider other isometric models as necessary: 
the {\it conformal ball model} $\BB^n$ and the {\it upper half-space model} $\UU^n$. 
For a subset $S\subset \BB^n$, 
$\overline{S}$ denotes the closure of $S$ in $\overline{\BB^n}$. 
A subset $S$ of $\BB^n$ is {\it convex} if, for each pair of
distinct points $x,y\in S$, the geodesic segment $[x,y]$ is contained in $S$.
For a convex open subset $S\subset\BB^n$ (or its closure in $\BB^n$), 
a {\it side} of $S$ is a maximal convex subset of the topological boundary $\partial S$. 
A {\it generalized polytope} in $\BB^n$
is a convex hull of finitely many points in $\overline{\BB^n}$. 
%
A {\it horoball} is an open ball in $\BB^n$ tangent to $\partial\BB^n$. 


Let ${\rm O}(L_\RR)$ (resp. ${\rm O}(L)$) be the orthogonal group of $L_\RR$ (resp. $L$). 
We write ${\rm O}^+(L_\RR)$ for the index two subgroup of transformations preserving the cone $\C$, 
which is naturally isomorphic to the isometry group $\Isom(\HH^n)$. 
We define a discrete subgroup ${\rm O}^+(L):={\rm O}^+(L_\RR)\cap{\rm O}(L)<\Isom(\HH^n)$. 
%
Each element $g\in\Isom(\HH^n)$ is classified as: 
\begin{itemize}
\item
{\it Elliptic}: $g$ fixes a point in $\HH^n$. 

\item 
{\it Parabolic}: $g$ fixes no point in $\HH^n$ and fixes a unique point in $\partial\HH^n$. 

\item 
{\it Loxodromic}: $g$ fixes no point in $\HH^n$ and fixes exactly two points in $\partial\HH^n$. 
\end{itemize}


Note that, for a fixed point 
of a parabolic (resp. loxodromic) isometry in ${\rm O}^+(L)$,
a representative fixed isotopic vector in $\partial\C$ is rational (resp. irrational) \cite[Remarque 1.1]{Cantat-K3}. 

\begin{dfn}
A subgroup of $\Isom(\HH^n)$ is {\it Kleinian} if it is discrete. 
\end{dfn}

A Kleinian group is {\it elementary} if it is virtually abelian, that is, it contains an abelian subgroup of finite index. 
Each elementary Kleinian group $\Gamma$ is also classified into the following three types: 
\begin{itemize}
\item
{\it Elliptic type}: $\Gamma$ fixes a point in $\HH^n$, or equivalently, $\Gamma$ is finite.  

\item 
{\it Parabolic type}: $\Gamma$ has a unique fixed point in $\partial\HH^n$, or equivalently, 
$\Gamma$ has a free abelian subgroup of finite index generated by parabolic elements. 

\item 
{\it Loxodromic type}: $\Gamma$ has two fixed point in $\partial\HH^n$, or equivalently,  
$\Gamma$ has an infinite cyclic subgroup of finite index generated by a loxodromic element. 
\end{itemize} 
Let $\Gamma$ be a Kleinian group. 
A point $c\in \partial\HH^n$ is a {\it limit point} of $\Gamma$ if  
there exists $x\in\HH^n$ and a sequence $\{g_i\}_{i=1}^\infty\subset\Gamma$
such that $\{g_ix\}_{i=1}^\infty$ converges to $c$. 
The {\it limit set} $\Lambda(\Gamma)$ of $\Gamma$ is the set of all limit points of $\Gamma$. 
Note that a fixed point of either a parabolic or loxodromic element of $\Gamma$ is a limit point of $\Gamma$. 
For each $x\in\HH^n$, 
we have $\Lambda(\Gamma)=\overline{\Gamma.x}\cap\partial\HH^n$, 
hence $\Lambda(\Gamma)$ is a $\Gamma$-invariant closed subset of $\partial\HH^n$. 

\begin{prop}
Let $\Gamma$ be a Kleinian group. 
The following are equivalent: 
\begin{enumerate}
\item
$\Gamma$ is elementary. 

\item
$\Lambda(\Gamma)$ is finite, especially $|\Lambda(\Gamma)|\le2$. 

\item 
$\Gamma$ has a finite orbit in $\overline{\HH^n}$. 
\end{enumerate}
\end{prop}

Thus, for each $c\in\partial\HH^n$, its stabilizer $\Gamma_c$ is elementary. 


The set $C(\Gamma)$ is the convex hull of $\Lambda(\Gamma)$ in $\HH^n$, that is, 
the smallest convex subset of $\HH^n$ satisfying $\overline{C(\Gamma)}\cap\partial\HH^n=\Lambda(\Gamma)$. 
%
$C(\Gamma)$ is $\Gamma$-invariant and closed in $\HH^n$, and 
when $\Gamma$ is non-elementary, 
any $\Gamma$-invariant closed subset of $\HH^n$ contains $C(\Gamma)$. 
The quotient $C(\Gamma)/\Gamma$ is called the {\it convex core} of the hyperbolic orbifold $\HH^n/\Gamma$. 



\vspace{3mm}
We introduce the notion of the geometrical finiteness, 
which is a key concept for our results. 

\begin{dfn}
\begin{enumerate}
\item
A Kleinian group $\Gamma<\Isom(\HH^n)$ is {\it geometrically finite} 
if $\Gamma$ is finitely generated, and 
$
{\rm vol}(C_\epsilon(\Gamma)/\Gamma)<\infty
$
for some $\epsilon>0$, where $C_\epsilon(\Gamma)$ is a $\epsilon$-neighborhood of $C(\Gamma)$ in $\HH^n$. 

\item
Let $G$ be a group. 
A representation $G\to\Isom(\HH^n)$ is {\it geometrically finite} if 
the image is a geometrically finite Kleinian group, and the kernel is finite. 
\end{enumerate}
\end{dfn}

The following is one of several characterizations of geometrically finite Kleinian groups. 

\begin{thm}[{\cite{Bow1}, see also \cite{Rat, Kap}}]\label{thm-trunc-conv-hull}
Let $\Gamma$ be a Kleinian group. 
The following are equivalent: 
\begin{enumerate}
\item 
$\Gamma$ is geometrically finite. 

\item 
There exists a $\Gamma$-invariant, pairwise disjoint collection $\{V_\lambda\}_\lambda$ of open horoballs at parabolic fixed points of $\Gamma$, such that the quotient 
$
\left(C(\Gamma)\backslash \bigcup_\lambda V_\lambda\right)/\Gamma
$
is compact. 
\end{enumerate}
\end{thm}


Let $\N\subset \HH^n$ be a convex subset
and $\Gamma$ a group of isometries of $\HH^n$ preserving $\N$. 
A connected closed subset $\Pi_\N$ of $\HH^n$ is a {\it fundamental domain} for the action of $\Gamma$ on $\N$ 
if the members of $\{g\Pi_\N^\circ~|~g\in\Gamma\}$ are mutually disjoint, and
$\N=\bigcup_{g\in\Gamma}g\Pi_\N$.
A fundamental domain $\Pi_\N$ for for the action of $\Gamma$ on $\N$ is {\it locally finite}
if $\{g\Pi_\N~|~g\in\Gamma\}$ is a locally finite collection of subsets of $\N$. 
A convex fundamental domain $\Pi_\N$ for the action of $\Gamma$ on $\N$ is {\it exact} if for each side $S$ of $\Pi_\N$ there is an element $g\in\Gamma$ such that $S=\Pi_\N\cap g\Pi_\N$.

\subsection{Groups and Cones in algebraic geometry}


Let $X$ be a smooth projective variety over a field $K$. 
The automorphism group (resp. birational automorphism group) of $X$ is denoted by $\Aut(X)$ (resp. $\Bir(X)$). 
We write $N^1(X)$ for the torsion-free part of the N\'{e}ron–Severi group $\NS(X)$. 
Its rank is called the Picard rank and denoted by $\rho_X$. 

We define several cones in $N^1(X)_\RR$ as follows:  
\begin{itemize}
\item
The {\it positive cone} $\C_X$ is the connected component of the set
$
\left\{v\in N^1(X)_\RR \mid v^2>0\right\}
$
containing ample classes. 
\item
The cone $\C^+_X\subset \overline{\C_X}$ is 
the convex hull of $\overline{\C_X}\cap N^1(X)_\QQ$ in $N^1(X)_\RR$. 
\item 
The {\it ample cone} $\Amp_X\subset \C_X$ is the cone generated by ample divisor classes.
\item 
The {\it nef cone} $\Nef_X\subset \overline{\C_X}$ is the closure of $\Amp_X$ in $N^1(X)_\RR$.
\item 
The cone $\Nef^+_X\subset \Nef_X$ is 
the convex hull of $\Nef_X\cap N^1(X)_\QQ$ in $N^1(X)_\RR$. 
\item 
The {\it effective cone} $\Eff_X\subset N^1(X)_\RR$ is the cone generated by integral curve classes. 
\item 
The {\it effective nef cone} $\Nef^e_X:=\Nef_X\cap\Eff_X$. 
\item 
The {\it movable cone} $\Mov_X\subset \overline{\C_X}$ is the cone generated by movable divisor classes.
\item 
The cone $\Mov^+_X\subset \overline{\Mov_X}$ is 
the convex hull of $\overline{\Mov_X} \cap N^1(X)_\QQ$ in $N^1(X)_\RR$. 
\end{itemize}
A closed cone $C\subset N^1(X)_\RR$ is {\it rational polyhedral} if
$C$ is the convex hull  of finitely many rational vectors, i.e. elements of $N^1(X)_\QQ$. 


Note that in any case of Theorem \ref{thm-geom-fin}, 
$N^1(X)$ admits a lattice structure of signature $(1,\rho_X-1)$: 
the intersection form for surfaces, and the Beauville--Bogomolov--Fujiki form for smooth irreducible symplectic varieties. 
%
%
For a rational polyhedral cone $C\subset N^1(X)_\RR$, 
the convex subspace $C\cap\HH^{\rho_X-1}$ is a generalized polytope whose vertices lie in $\overline{\HH^{\rho_X-1}}(\QQ)$. 


\section{Proof of Theorem \ref{thm-geom-fin}}\label{pf-main-thm}
\subsection{General statement}
\label{sec-gen-statement}
We shall prove the following general statement on the geometrical finiteness. 

Throughout this subsection, let $L$ be a lattice of signature $(1,n)$ and 
fix a positive cone $\C\subset L_\RR$. 

\begin{thm}\label{thm-general-statement}
Let $\Gamma$ be a subgroup of ${\rm O}^+(L)$ 
and 
$\N$ be a $\Gamma$-invariant closed convex subset of $\HH^n$. 
Suppose that there exists 
a fundamental domain $\Pi_\N$ for the action of $\Gamma$ on $\N$
satisfying the following conditions: 
\begin{enumerate}
\item 
$\Pi_\N$ is locally finite and exact, 

\item 
$\Pi_\N$ is a generalized polytope whose vertices lie in $\overline{\HH^n}(\QQ)$. 

\end{enumerate}
Then $\Gamma$ is geometrically finite. 
\end{thm}

To prove this theorem, we take a similar approach as in \cite[Section 12.4]{Rat}. 
Elementary Kleinian groups are geometrically finite. 
Therefore, in the following, 
we assume that $\Gamma$ is non-elementary. 


Note that 
if the set $\overline{\Pi}_\N\cap \Lambda(\Gamma)$ is nonempty, then it is finite since $\Pi_\N$ is a generalized polytope. 
The following is a key to prove Theorem \ref{thm-general-statement}. 

\begin{lem}
If $\overline{\Pi}_\N\cap \Lambda(\Gamma)$ is nonempty, then it consists of parabolic fixed points of $\Gamma$. 
\end{lem}
\begin{pf}
Since any $c\in\overline{\Pi}_\N\cap \Lambda(\Gamma)$ is rational, 
the stabilizer $\Gamma_c$ is elementary of either elliptic or parabolic type, 
hence it is enough to show that $\Gamma_c$ is infinite. 


We pass to the upper half-space model $\UU^n$ and conjugate $\Gamma$ so that $c=\infty$. 
Let ${\rm v}:\UU^n\to\EE^{n-1}$ be the vertical projection. 
We define a subset $U$ of $\N$ as follows 
\[
U:=\cup\{g\Pi_\N \mid g\in\Gamma\text{ such that }c\in g\overline{\Pi}_\N\}. 
\]
We now show that ${\rm v}U={\rm v}\N$. 
Since $\{g\Pi_\N\}_{g\in\Gamma}$ is locally finite, ${\rm v}U$ is closed in ${\rm v}\N$. 
Let us show that ${\rm v}U$ is open in ${\rm v}\N$. 
For any $z\in{\rm v}U$, we take an element $w\in{\rm v}^{-1}z$. 
If $w$ is an inner point of $U$, then $z$ is also. 
If not, 
we may assume that 
$w\in\partial U$ lies in a vertical side (i.e. a side whose closure contains $c=\infty$) of $f_1\Pi_\N$ for some $f_1\in\Gamma$. 
There exist $f_2,\cdots,f_s\in \Gamma$ and 
a sufficiently small open ball $B(w)$ at $w$ 
such that $w\in f_i\Pi_\N$ for each $i$ and 
\[
\N\cap B(w)\subset \cup^s_{j=1} f_j\Pi_\N. 
\]
Furthermore, 
we can assume that
$B(w)$ does not meet non-vertical sides of each $f_i\Pi_\N$. 
Note that \cite[Theorem 6.7.5]{Rat} is available since $\Pi_\N$ is a locally finite, exact fundamental domain and a generalized polytope. 
For some $f_{i_1}\in\Gamma$, 
if $f_{i_1}\Pi_\N$ has a side that intersects a vertical one of $f_1\Pi_\N$ containing $w$, 
then 
the side of $f_{i_1}\Pi_\N$ is also vertical 
by \cite[Theorem 6.7.5]{Rat}. 
Hence, any side of $f_{i_1}\Pi_\N$ containing $w$ is vertical. 
If some $f_{i_2}\Pi_\N$ has a side that intersects a vertical one of $f_{i_1}\Pi_\N$ containing $w$, 
the same argument implies that any side of $f_{i_2}\Pi_\N$ containing $w$ is vertical. 
Repeating this argument, we find that
each $f_i\overline{\Pi}_\N$ contains $c$, 
namely $\cup^s_{j=1} f_j\Pi_\N\subset U$. 
We have 
\[
z\in 
{\rm v}\N\cap{\rm v}B(w)
={\rm v}(\N\cap B(w))
\subset \cup^s_{j=1} {\rm v}f_j\Pi_\N
\subset {\rm v}U.
\]
It turns out that ${\rm v}U$ is open in ${\rm v}\N$. 
Hence ${\rm v}U={\rm v}\N$. 


By the finiteness of $\overline{\Pi}_\N\cap \Gamma.c$, 
we have elements $h_1,\cdots,h_k\in \Gamma$ such that
$\overline{\Pi}_\N\cap \Gamma.c=\{h^{-1}_jc\mid j=1,\cdots,k\}$. 
Note that, for $g\in\Gamma$,  
we have $c\in g\overline{\Pi}_\N$ if and only if $g\in\cup^k_{j=1}\Gamma_c h_j$. 
Therefore, 
\[
{\rm v}\N=\cup\left\{{\rm v}g\Pi_\N~\middle|~g\in\cup^k_{j=1}\Gamma_c h_j\right\}. 
\]


Suppose to the contrary that $\Gamma_c$ is finite. 
Then we have $\cup^k_{j=1}\Gamma_c h_j=\{g_1,\cdots,g_m\}$ for some elements of $\Gamma$. 
For any $y\in\N$, 
there exists $y_0\in g_{j_0}\Pi_\N$ for some $g_{j_0}\in \Gamma$
such that 
the loxodromic geodesic ray, or equivalently, the vertical ray $[y,c)\subset\UU^n$ satisfies that
\[
[y,c)\cap g_{j_0}\Pi_\N=[y_0,c).
\]
Thus, taking a sufficiently small open ball $B(c)$ at $c$ in $\overline{\UU^n}$, 
we have $\N\cap B(c)\subset\cup^m_{j=1}g_j\Pi_\N$. 
By $c\in \Lambda(\Gamma)=\overline{\Gamma.x}\cap\partial\HH^n$
for any $x\in g_1\Pi^\circ_\N$, 
there exist a infinite sequence $\{f_j\}_{j\ge1}\subset\Gamma$ 
such that
the sequence $\{f_jx\}_j\subset \N\cap B(c)$ converges to $c$, 
hence $\{f_jx\}_j\subset\cup^m_{j=1}g_j\Pi_\N$. 
However, this contradicts that $\Pi_\N$ is a fundamental domain. 
\qed
\end{pf}

Set $\overline{\Pi}_\N\cap \Lambda(\Gamma)=\{c_1,\cdots,c_l\}\subset\partial\HH^n(\QQ)$. 

\begin{lem}[{\cite{Tot}}]
There exists an open horoball $B_i$ 
at each $c_i\in\overline{\Pi}_\N\cap \Lambda(\Gamma)$ 
such that 
elements in the collection
\[\left\{gB_i \mid g\in\Gamma,~i=1,\cdots,l\right\}\]
are pairwise disjoint. 
\end{lem}
\begin{pf}
The claim is proved in \cite[Section 2]{Tot}, but we provide a proof here for reader convenience. 
For $c_i\in \overline{\Pi}_\N\cap \Lambda(\Gamma)$, 
let $e_i\in \partial\C\cap L$ be the corresponding primitive isotropic integral vector. 
We explicitly define an open horoball $B_i$ at $c_i$ as follows
\[
B_i:=B_{e_i}:=\left\{x\in\HH^n ~\middle|~ (x,e_i)<\frac{1}{2}\right\}. 
\]
For distinct primitive isotropic integral vectors $e,e'\in \partial\C\cap L$ and $x\in\HH^n$, by \cite[Theorem 3.1.1]{Rat} and \cite[Lemma 2.1]{Tot}, we have
\[
1\le(e,e')\le 2(x,e)(x,e'),
\]
which implies that $B_e$ and $B_{e'}$ are disjoint. 
Furthermore, it is easy to check that $gB_{e}=B_{ge}$ for each $g\in\Gamma$, 
which completes the proof. 
\qed
\end{pf}

Therefore, we obtain
the $\Gamma$-invariant, pairwise disjoint collection
$\B:=\left\{gB_i \mid g\in\Gamma,~i=1,\cdots,l\right\}$ 
of open horoballs at parabolic fixed points of $\Gamma$.

\vspace{3mm}
\noindent
{\it Proof of Theorem \ref{thm-general-statement}.}
\vspace{2mm}

\noindent
Let $\pi:\HH^n\to\HH^n/\Gamma$ be the quotient map. 
Then it easily follows from $C(\Gamma) \subset \N$ that
\begin{equation}\label{eq-pf-gen-stat}
\pi\left(C(\Gamma)\backslash \bigcup_\B gB_i \right)
\subset
\pi\left(\Pi_\N \cap C(\Gamma)\backslash \bigcup_{i=1}^l B_i\right)
. 
\end{equation}
Assume that $\Pi_\N \cap C(\Gamma)\backslash \bigcup_{i=1}^l B_i$  is unbounded. 
Then there exists unbounded sequence $\{a_j\}_{j\ge1} \subset \Pi_\N \cap C(\Gamma)\backslash \bigcup_{i=1}^l B_i$. 
We can take a subsequence $\{a'_j\}$
convergent to a point $y$ in the set
\[
\left(\overline{\Pi_\N \cap C(\Gamma)}\right)\cap\partial\HH^n
\subset
\overline{\Pi}_\N\cap \Lambda(\Gamma)
=\{c_1,\cdots,c_l\}, 
\]
hence $y=c_{i_0}$ for some $i_0$. 
Since $\Pi_\N$ is finite-sided, 
$a'_j$ is in the horoball $B_{i_0}$ for sufficiently large $j$, 
which contradicts the assumption. 
Therefore, the right-hand side of (\ref{eq-pf-gen-stat}) is compact, 
thus the left-hand side is also. 
By Theorem \ref{thm-trunc-conv-hull}, $\Gamma$ is geometrically finite.  
\qed

\subsection{Cone conjecture implies Geometrical finiteness}

Applying Theorem \ref{thm-general-statement} and the cone conjecture, 
we now complete the proof of Theorem \ref{thm-geom-fin} for each case. 


\vspace{3mm}
\noindent
{\it Proof of Theorem \ref{thm-geom-fin}.}

\vspace{2mm}
\noindent
For simplicity, let us denote $\tau:=\rho_X-1$. 

\vspace{2mm}
\noindent
{\bf Case (1): }
Let $X$ be a K3 surface over a field of characteristic different from $2$. 
Let us recall that the representation $\Aut(X)\to{\rm O}^+(\NS(X))$ 
has a finite kernel (see \cite{HuyK3, Lieblich-Maulik, BLvL}), and the image is denoted by $\Gamma_X$.
For the Weyl group $W_X$ of $X$, 
we define $\widetilde{\Gamma}_X:=\Gamma_X\ltimes W_X$. 

We apply Theorem \ref{thm-general-statement} with 
$L:=\NS(X)$ and $\N:=\Nef_X^+\cap\HH^\tau$. 
Fix an element $h:=\frac{H}{\sqrt{(H,H)}}\in\N$ for some ample class $H\in\NS(X)$ such that the stabilizer $\Gamma_{X,H}$ is trivial. 
We here recall the cone conjecture: 
\begin{thm}[{\cite{Sterk, TotDuke, Lieblich-Maulik, BLvL, Takamatsu}}]
\label{thm-cone-conj-K3}
Let $X$ be a K3 surface over a field $K$. 
Suppose that the characteristic of $K$ is different from $2$ or $X$ is not supersingular. 
Then
\[
D^+_X:=\left\{x\in \C^+_X ~\middle|~(H,x)\le(H,gx)\text{ for any }g\in \widetilde{\Gamma}_X\right\}
\]
is a rational polyhedral fundamental domain for the action of $\Gamma_X$ on $\Nef^+_X$. 
\end{thm}
A closed subset $\Pi'_\N:=D^+_X\cap\HH^\tau$
is a fundamental domain for the action of $\Gamma_X$ on $\N$. 
Clearly, $\Pi'_\N$ is a subset of the Dirichlet domain (\cite[p240]{Rat}),  
which is a locally finite fundamental domain for the action of $\Gamma_X$ 
on $\HH^\tau$. 
Thus, $\Pi'_\N$ is also locally finite 
which implies that $\N$ is closed in $\HH^\tau$, hence is proper, i.e. any closed ball is compact. 
$\N$ is also geodesically connected and geodesically complete in the sense of \cite{Rat}. 
%
By \cite[Theorem 6.6.13 and Theorem 6.7.4]{Rat}, 
the Dirichlet domain $\Pi_\N$ 
defined by
\[
\Pi_\N:=
\left\{x\in \N ~\middle|~ d(h,x)\le d(h,gx)\text{ for any }g\in \Gamma_X\right\},
\]
is a locally finite, exact fundamental domain for the action of $\Gamma_X$ on $\N$. 
Therefore, the inclusion $\Pi'_\N\subset\Pi_\N$ is actually an equality. 
%
Since $\Pi'_\N$ is also a generalized polytope whose vertices lie in 
$\overline{\HH^\tau}(\QQ)$, 
$\Pi_\N$ satisfies all the conditions in Theorem \ref{thm-general-statement}. 


\vspace{3mm}
\noindent
{\bf Case (2): }
The proof is the same as Case (1).

\vspace{3mm}
\noindent
{\bf Case (3): }
Let $X$ be an Enriques surface over an algebraically closed field of characteristic different from $2$. 
The representation $\Aut(X)\to{\rm O}^+(N^1(X))$ 
has a finite kernel (see \cite[Proposition 2.1]{DolgachevMartin2019}), and the image is denoted by $\Gamma_X$.
For the subgroup $W^{{\rm nod}}_X$ of ${\rm O}^+(N^1(X))$ generated by reflections associated with classes of nodal curves in $X$, 
we define $\widetilde{\Gamma}_X:=\Gamma_X\ltimes W^{{\rm nod}}_X$. 

The proof is similar to Case (1).
We apply Theorem \ref{thm-general-statement} with 
$L:=N^1(X)$ and $\N:=\Nef_X^e\cap\HH^9$. 
The cone conjecture is as follows: 
\begin{thm}[{\cite{Namikawa, Kawamata, Oguiso-Sakurai, Wang}}]
Let $X$ be an Enriques surface over an algebraically closed field of characteristic different from $2$. 
Then
\[
D^+_X:=\left\{x\in \C^+_X ~\middle|~(H,x)\le(H,gx)\text{ for any }g\in \widetilde{\Gamma}_X\right\}
\]
is a rational polyhedral fundamental domain for the action of $\Gamma_X$ on $\Nef^e_X$. 
\end{thm}

\vspace{3mm}
\noindent
{\bf Case (4): }
In this paper, 
a Coble surface means a {\it terminal Coble surface of K3 type} in the sense of Dolgachev--Zhang \cite{Dolgachev-Zhang}. 
Let $X$ be a Coble surface over an algebraically closed field of characteristic $0$, and 
$Y\to X$ the K3 cover with covering involution $\sigma$. 
A natural representation 
\[
\Aut(X)\to{\rm O}^+(\Pic(X))\simeq{\rm O}^+(\NS(Y)^\sigma)<{\rm O}^+(\NS(Y)) 
    \]
has a finite kernel (see \cite{Cantat-Dolgachev, FulinXu}), 
and the image in ${\rm O}^+(\NS(Y)^\sigma)$ is denoted by $\Gamma_X$.
For the subgroup $R_X$ of ${\rm O}^+(\NS(Y)^\sigma)$
called $\sigma$-equivariant reflection group (\cite{Oguiso-Sakurai}), 
we define $\widetilde{\Gamma}_X:=\Gamma_X\ltimes R_X$. 

The proof is similar to Case (1).
We apply Theorem \ref{thm-general-statement} with $L:=\NS(Y)^\sigma\simeq\Pic(X)$, 
$\C_X:=\C_Y\cap\NS(Y)^\sigma$, 
and $\N:=(\Nef_Y^e)^\sigma\cap\HH^\tau \simeq \Nef_X^e\cap\HH^\tau$. 
The cone conjecture is as follows: 
\begin{thm}[{\cite{Oguiso-Sakurai, FulinXu}}]
Let $X$ be a Coble surface 
over an algebraically closed field of characteristic $0$. 
Then
\[
D^+_X:=\left\{x\in \C^+_X ~\middle|~(H,x)\le(H,gx)\text{ for any }g\in \widetilde{\Gamma}_X\right\}
\]
is a rational polyhedral fundamental domain for the action of $\Gamma_X$ on $(\Nef^e_Y)^\sigma$. 
\end{thm}

\vspace{3mm}
\noindent
{\bf Case (5): }
Let $X$ be a smooth irreducible symplectic variety over a field of characteristic $0$. 
The representation $\Aut(X)\to{\rm O}^+(N^1(X))$ 
has a finite kernel (see \cite{Takamatsu}), 
and the image is denoted by $\Gamma_X$. 

We apply Theorem \ref{thm-general-statement} with 
$L:=N^1(X)$ and $\N:=\Nef_X^+\cap\HH^\tau$. 
Fix an element $h:=\frac{H}{\sqrt{(H,H)}}\in\N$ for some ample class $H\in N^1(X)$ such that the stabilizer $\Gamma_{X,H}$ is trivial. 
The cone conjecture is as follows: 
\begin{thm}[{\cite{Amerik-Verbitsky, Takamatsu}}]
Let $X$ be a smooth irreducible symplectic variety over a field of characteristic $0$. 
Then
\[
D^+_X:=\left\{x\in \Nef_X ~\middle|~(H,x)\le(H,gx)\text{ for any }g\in\Gamma_X\right\}
\]
is a rational polyhedral fundamental domain for the action of $\Gamma_X$ on $\Nef^+_X$. 
\end{thm}
In this case, 
a closed subset $\Pi_\N:=D^+_X\cap\HH^\tau$
is nothing but the Dirichlet fundamental domain for the action of $\Gamma_X$ on $\N$. 
The remaining arguments are the same as Case (1). 

\vspace{3mm}
\noindent
{\bf Case (6): }
Let $X$ be a smooth irreducible symplectic variety over a field $K$ of characteristic $0$. 
The representation $\Bir(X)\to{\rm O}^+(N^1(X))$ 
has a finite kernel (see \cite{Takamatsu}), 
and the image is denoted by $\Gamma_X$.
Let $W^{{\rm Exc}}_{X_{\overline{K}}}$ be the subgroup of ${\rm O}^+(N^1(X_{\overline{K}}))$ generated by reflections associated with classes of prime exceptional divisors on $X_{\overline{K}}$, and 
$R_X$ the ${\rm Gal}_K$-fixed part in $W^{{\rm Exc}}_{X_{\overline{K}}}$
(\cite{Takamatsu}), 
where ${\rm Gal}_K$ is the absolute Galois group of $K$. 
Note that $R_X$ faithfully acts on $N^1(X)$. 
We define $\widetilde{\Gamma}_X:=\Gamma_X\ltimes R_X<{\rm O}^+(N^1(X))$. 

The proof is similar to Case (1).
We apply Theorem \ref{thm-general-statement} with 
$L:=N^1(X)$ and $\N:=\Mov_X^+\cap\HH^\tau$. 
The cone conjecture is as follows: 
\begin{thm}[{\cite{Markman, Takamatsu}}]
Let $X$ be a smooth irreducible symplectic variety over a field of characteristic $0$. 
Then
\[
D^+_X:=\left\{x\in \C^+_X ~\middle|~(H,x)\le(H,gx)\text{ for any }g\in \widetilde{\Gamma}_X\right\}
\]
is a rational polyhedral fundamental domain for the action of $\Gamma_X$ on $\Mov_X^+$. 
\end{thm}
\qed

\section{Applications}

We present several applications of Theorem \ref{thm-geom-fin}. 

\subsection{Non-positively curved properties}

Throughout this subsection, let $G_X$ be a group as in Theorem \ref{thm-geom-fin}. 
We first observe that $G_X$ is non-positively curved: relatively hyperbolic and ${\rm CAT(0)}$, see \cite{Bow2, Osin, BriH} for definitions. 



Historically, non-elementary geometrically finite Kleinian groups have served as the original examples of relatively hyperbolic groups with respect to their actions on the convex hulls of their limit sets, where peripheral subgroups are given by maximal parabolic subgroups (\cite{Bow2}). 
Recall that 
a (relatively) hyperbolic group is called {\it elementary} if it is virtually cyclic. 

\begin{cor}[Alternative]\label{cor-alternative}
$G_X$ is either virtually abelian or non-elementary relatively hyperbolic. 
\end{cor}
\begin{pf}
The claim is clear by Theorem \ref{thm-geom-fin}. 
Note that since the isometric action on $\HH^{\rho_X-1}$ has a finite kernel, 
$\Gamma_X$ is non-elementary relatively hyperbolic (resp. virtually abelian) 
if and only if 
$G_X$ is non-elementary relatively hyperbolic (resp. virtually abelian). 
\end{pf}

We also obtain a criterion for relative hyperbolicity. 

\begin{cor}
$G_X$ is non-elementary relatively hyperbolic if and only if $G_X$ contains the rank $2$ free group $F_2$. 
\end{cor}
\begin{pf}
This is a direct corollary of the strong Tits alternative: 
for every subgroup $H<G_X$, either $H$ is virtually abelian or $H$ contains $F_2$. 
Note that relatively hyperbolic groups satisfy the Tits alternative (\cite{Tuk}), and 
the virtual solvability is equivalent to the virtual abelianity 
in this case by \cite[Ch.I\hspace{-1.2pt}I, Theorem 7.8]{BriH}. 
\qed
\end{pf}


\begin{cor}\label{cor-cat(0)}
$G_X$ is {\rm CAT(0)}. 
\end{cor}
\begin{pf}
It is well-known that each geometrically finite Kleinian group $\Gamma$ is ${\rm CAT(0)}$ 
via the isometric action on the ${\rm CAT(0)}$ space $C(\Gamma)\backslash \bigcup_\lambda V_\lambda$
with the induced length metric,  
called the {\it truncated convex hull} of the limit sets, 
where $\{V_\lambda\}$ is a $\Gamma$-invariant pairwise disjoint collection of open horoballs at parabolic fixed points of $\Gamma$ (\cite[Ch.I\hspace{-1.2pt}I, Exercise 11.37(2)]{BriH}). 

Since the representation $G_X\to\Isom(\HH^{\rho_X-1})$ has a finite kernel, 
$G_X$ is also ${\rm CAT(0)}$ through its proper and cocompact isometric action on the truncated convex hull.
\qed
\end{pf}

Especially, each $G_X$ is finitely presented by \cite[Ch.I\hspace{-1.2pt}I\hspace{-1.2pt}I.$\Gamma$, Theorem 1.1]{BriH}. 

\subsection{Dynamical characterizations for K3 surfaces}

In this subsection, let $X$ be a K3 surface over an algebraically closed field $K$. 
We further suppose that the characteristic of $K$ is different from $2$, or $X$ is not supersingular as in Theorem \ref{thm-geom-fin}. 

Entropy provides a dynamical perspective on the complexity and behavior of automorphism groups. 
For each automorphism $f\in\Aut(X)$, 
the {\it entropy} of $f$ 
is numerically defined as the logarithm of the spectral radius of its induced action on $\NS(X)_\RR$. 
Then an automorphism of infinite order is parabolic (resp. loxodromic) 
if and only if 
its entropy is zero (resp. positive). 
%
We call that
{\it $\Aut(X)$ has zero entropy} if the entropy of any automorphism of $X$ is zero, 
otherwise, {\it $\Aut(X)$ has positive entropy}. 
Note that, over the complex number field, 
the entropy is actually equal to the {\it topological entropy} due to Gromov--Yomdin \cite{Gro,Yom}. 

The following is a characterization of the virtual abelianity via entropy. 

\begin{thm}[{\cite{Cantat-K3, Oguiso-20-zero, Brandhorst-Mezzedimi, Yu}}]\label{zero-entropy-classification}
Suppose that $\Aut(X)$ is infinite. 

\begin{enumerate}
\item
For $\rho_X=2$, 
$\Aut(X)$ is virtually cyclic and has positive entropy. 

\item 
For $\rho_X=3$, 
$\Aut(X)$ is virtually abelian if and only if
$\Aut(X)$ is virtually cyclic. 

\item 
For $\rho_X=4$, 
$\Aut(X)$ is virtually abelian if and only if
either $\Aut(X)$ is virtually cyclic and has positive entropy, 
or $\Aut(X)$ has zero entropy. 

\item 
For $\rho_X\ge 5$, 
$\Aut(X)$ is virtually abelian if and only if 
$\Aut(X)$ has zero entropy. 

\end{enumerate}
\end{thm}
\begin{pf}
See \cite[Corollary 1.5, Remark 3.10, Proposition 3.11]{Brandhorst-Mezzedimi}, 
but note that we use the hyperbolicity (see Example \ref{ex-hyp}) in the case of $\rho_X=2,3$. 
\qed
\end{pf}

By the alternative in Corollary \ref{cor-alternative}, 
we obtain a dynamical and numerical characterization of relative hyperbolicity. 

\begin{cor}\label{cor-dyn}
Suppose that the Picard rank of $X$ is at least $5$. 
Then $\Aut(X)$ is non-elementary relatively hyperbolic if and only if 
$\Aut(X)$ has positive entropy. 
\end{cor}

\begin{rem}
A statement similar to Corollary \ref{cor-dyn} holds for Enriques surfaces by combining Corollary \ref{cor-alternative} with a 
result by Martin--Mezzedimi--Veniani \cite[Proposition 2.8]{MMV}.
\end{rem}

\subsection{Lattices and Convex-cocompact groups for K3 surfaces}

Lattices and convex-cocompact groups 
form important classes within geometrically finite Kleinian groups.  
For further details, we refer to \cite{Kap} and the references therein. 
Two theorems (Theorem \ref{chara-lattice}, \ref{chara-conv-cocpt}) below claim that these notions are related to $(-2)$-curves and genus-one fibrations respectively. 

Let $X$ be a K3 surface over an algebraically closed field $K$. 
We further suppose that the characteristic of $K$ is different from $2$, or $X$ is not supersingular as in Theorem \ref{thm-geom-fin}, 
and that $\Aut(X)$ is non-elementary, hence $\rho_X$ is at least $3$.  

\begin{dfn}
\begin{enumerate}
\item
A Kleinian group $\Gamma<\Isom(\HH^n)$ is a {\it lattice} 
if 
$\Lambda(\Gamma)=\partial\HH^n$ and $\Gamma$ is geometrically finite. 
\item
Let $G$ be a group. 
A representation $G\to\Isom(\HH^n)$ is a {\it lattice} if 
the image is a lattice, and the kernel is finite. 
\end{enumerate}
\end{dfn}

\begin{thm}\label{chara-lattice}
The representation $\Aut(X)\to\Isom(\HH^{\rho_X-1})$ is a lattice 
if and only if $X$ has no $(-2)$-curves. 
\end{thm}
\begin{pf}
Note that $X$ has a $(-2)$-curve if and only if $\NS(X)$ has a root, 
and we have
\[
\Amp_X=\left\{v\in\C_X \mid (v,C)>0\text{ for all }(-2)\text{ curves }C\subset X \right\}, 
\]
see \cite[Ch.8, Corollary 1.7]{HuyK3}. 
If $\NS(X)$ has a root, then we have
\[
\Lambda(\Aut(X))
\subset \overline{\Amp_X\cap~ \HH^{\rho_X-1}}\cap \partial \HH^{\rho_X-1}
\subsetneqq \partial \HH^{\rho_X-1}
\]
since the complement of the open subset $\Amp_X\cap~ \HH^{\rho_X-1}$ in $\HH^{\rho_X-1}$ is unbounded. 
Suppose that $\NS(X)$ has no roots. 
Then the convex subset $\N$ in the proof of Theorem \ref{thm-geom-fin} is equal to $\HH^{\rho_X-1}$. 
By the cone conjecture (see Theorem \ref{thm-cone-conj-K3}), 
$\HH^{\rho_X-1}$ admits a fundamental domain $\Pi$ that is a generalized polytope. 
Therefore, we have 
\[
{\rm vol}(\HH^{\rho_X-1}/\Aut(X))={\rm vol}(\Pi)<\infty,
\]
hence $\Lambda(\Aut(X))=\partial\HH^{\rho_X-1}$ by \cite[Theorem 12.4.10]{Rat}. 
\qed
\end{pf}


\begin{dfn}
\begin{enumerate}
\item
A Kleinian group $\Gamma<\Isom(\HH^n)$ is {\it convex-cocompact} 
if 
there exists a nonempty closed convex $\Gamma$-invariant subset $C\subset\HH^n$ 
such that the quotient $C/\Gamma$ is compact. 
\item
Let $G$ be a group. 
A representation $G\to\Isom(\HH^n)$ is {\it convex-cocompact} if 
the image is a convex-cocompact Kleinian group, and the kernel is finite. 
\end{enumerate}
\end{dfn}

Note that convex-cocompact Kleinian groups are hyperbolic. 
Convex-cocompact lattices are called {\it uniform lattices}. 

\begin{prop}
Let $\Gamma$ be a Kleinian group. 
The following are equivalent: 
\begin{enumerate}
\item
$\Gamma$ is convex-cocompact. 
\item
The quotient $C(\Gamma)/\Gamma$ is compact. 
\item 
$\Gamma$ is geometrially finite with no parabolic elements. 

\end{enumerate}
\end{prop}

We remark that an infinite automorphism of $X$ is parabolic 
if and only if 
it preserves a genus-one fibration (i.e. stabilizes the fiber class of an elliptic or quasi-elliptic fibration) on $X$ (\cite{Cantat-K3, Oguiso-20-zero, Brandhorst-Mezzedimi}). 
The above characterization of convex-cocompactness implies the following. 

\begin{thm}\label{chara-conv-cocpt}
The representation $\Aut(X)\to\Isom(\HH^{\rho_X-1})$ is convex-cocompact 
if and only if 
either
\begin{enumerate}
\item
$X$ admits no genus-one fibrations, or
\item 
$X$ admits a genus-one fibration, and the Mordell--Weil group of the Jacobian fibration of any genus-one fibration is finite. 
\end{enumerate}
\end{thm}
\begin{pf}
The claim directly follows from Theorem \ref{thm-geom-fin}. 
Note that, for a genus-one fibration $\varphi: X\to \PP^1$,  
the Mordell--Weil group of the Jacobian fibration of $\varphi$ 
can be embedded into the stabilizer of $\varphi$ (i.e. a maximal parabolic subgroup)
as a finite index subgroup (see, e.g. \cite[Theorem 3.3]{Dolgachev-Martin} or \cite[Proposition 3.2]{Brandhorst-Mezzedimi}). 
\qed
\end{pf}

\begin{cor}\label{K3-conv-cocpt}
If $\rho_X\ge6$, 
then the representation $\Aut(X)\to\Isom(\HH^{\rho_X-1})$ is not convex-cocompact. 
Furthermore, if $\rho_X=5$, 
then this representation is convex-cocompact 
if and only if 
$\NS(X)$ is isomorphic to one of the following lattices:
\[
\langle 2^k\rangle \oplus D_4,
\qquad
\langle 2\cdot 3^{2m-1}\rangle \oplus A_2^{\oplus 2},
\]
where $k\ge 5$ and $m\ge 2$.  
\end{cor}
\begin{pf}
This follows immediately from Theorem \ref{chara-conv-cocpt} and 
Nikulins's classification of a N\'{e}ron--Severi lattice whose each Mordell--Weil group is finite (\cite[Theorem 10.2.2 a)]{Nikulin83}).
\qed
\end{pf}

\begin{cor}\label{K3-uniform-lattice}
\begin{enumerate}
\item 
If $\rho_X\ge5$, 
then the representation $\Aut(X)\to\Isom(\HH^{\rho_{X}-1})$ is not a uniform lattice. 

\item
There exists a complex K3 surface $X'$ of Picard rank $3$ or $4$
such that
the representation $\Aut(X')\to\Isom(\HH^{\rho_{X'}-1})$ is a uniform lattice. 
Furthermore, $\Aut(X')$ is virtually a fundamental group of a closed hyperbolic manifold of dimension $2$ or $3$, respectively. 


\end{enumerate}
\end{cor}
\begin{pf}
The claim (i) follows since each of lattices $\langle 2^k\rangle \oplus D_4,~
\langle 2\cdot 3^{2m-1}\rangle \oplus A_2^{\oplus 2}~(k\ge 5,~m\ge 2)$
in Corollary \ref{K3-conv-cocpt} has a root. 
There exists a complex K3 surface $X'$ with $\rho_{X'}=3$ or $4$
whose N\'{e}ron--Severi lattice contains no roots and no nonzero isotropic vectors. 
Then the representation $\Aut(X')\to\Isom(\HH^{\rho_{X'}-1})$ is a uniform lattice by Theorems \ref{chara-lattice} and \ref{chara-conv-cocpt}. 
This representation has a section whose image is the symplectic automorphism group $\Aut_s(X')$. 
Then $\Aut_s(X')\subset \Isom(\HH^{\rho_{X'}-1})$ is a matrix group, thus has a torsion-free subgroup $\Gamma$ of finite index. 
Since $\Gamma$ is a torsion-free uniform lattice, we clearly have $\Gamma\simeq\pi_1(\HH^{\rho_{X'}-1}/\Gamma)$, which completes the proof of the claim (ii). 
\qed
\end{pf}

In Example \ref{eg-vir-hyp-mfd-gp}, we see explicit examples of such lattices as in Corollary \ref{K3-uniform-lattice} (ii). 

\begin{rem}
In contrast to the case of K3 surfaces, the automorphism group of an Enriques surface is finite 
if and only if 
every Mordell--Weil group is finite \cite[Proposition 2.7, Remark 2.9]{MMV}. 
Consequently, the infinite automorphism group of an Enriques surface is never convex-cocompact. 
\end{rem}

\subsection{Examples}

We collect several examples in this subsection. 
\begin{eg}[{\cite{Brandhorst-Mezzedimi, Yu}}]
If the automorphism group $\Aut(X)$ of a K3 surface $X$ is virtually abelian, then $\rho_X\leq 18$, and the rank of $\Aut(X)$ is at most $8$. 
Conversely, for each integer $\rho'\in\{3,4,\cdots,18\}$, 
there exists a complex K3 surface of Picard rank $\rho'$
whose automorphism group is virtually abelian. 
\end{eg}


\begin{eg}\label{ex-hyp}
Let $X$ be a K3 surface as in Theorem \ref{thm-geom-fin}. 
If its Picard rank is at most $3$, 
then $\Aut(X)$ is (possibly elementary) hyperbolic by \cite[Ch.I\hspace{-1.2pt}I\hspace{-1.2pt}I.$\Gamma$, Theorem 3.1]{BriH}. 
\end{eg}


\begin{eg}[{\cite{Oguiso-20-zero, Yu-supersingular, Brandhorst-supersingular}, \cite[Corollary 5.4]{Brandhorst-Mezzedimi} and \cite[Theorem 0.3]{Mezzedimi-19}}]
The following K3 surfaces over an algebraically closed field have positive entropy, 
and hence their automorphism groups are non-elementary relatively hyperbolic. 
\begin{enumerate}
\item
Kummer surfaces in characteristic different from $2$. 
\item
K3 surfaces covering an Enriques surface, 
unless $\NS(X)\simeq U\oplus E_8\oplus D_8$
in characteristic different from $2$. 
\item 
complex K3 surfaces of Picard rank 19 with infinite automorphism group. 

\item 
Singular K3 surfaces.
\item 
Supersingular K3 surfaces in characteristic different from $2$. 
\end{enumerate}
\end{eg}


\begin{eg}
We shall consider two specific examples of singular K3 surfaces over $\CC$. 
\begin{enumerate}
\item
Let $X_F$ be the Fermat quartic. 
Then by Shioda \cite[Proposition 19]{Shi}, 
$X_F$ admits an elliptic fibration of the Mordell--Weil rank 6. 
Therefore, 
$\Aut(X_F)$ is relatively hyperbolic, but not hyperbolic. 

\item 
Let $X_3$ and $X_4$ be the K3 surfaces
whose transcendental lattices are of the form 
\[
\Trans(X_3)
=
\begin{pmatrix}
2 & 1 \\
1 & 2 \\
\end{pmatrix}
\text{  and  }
\Trans(X_4)
=
\begin{pmatrix}
2 & 0 \\
0 & 2 \\
\end{pmatrix}
\]
respectively, see \cite{SI} and \cite[Corollary 14.3.21]{HuyK3}. 
By Vinberg (\cite[Theorem in 2.4 and Theorem in 3.3]{Vin}), 
their automorphism groups are virtually free, hence hyperbolic. 
Furthermore, singular K3 surfaces have a parabolic element induced by a section of an elliptic fibration. 
Thus, $\Aut(X_3)$ (resp. $\Aut(X_4)$) has both distinct structures of hyperbolicity and relative hyperbolicity. 
In other words, 
$\Aut(X_3)$ (resp. $\Aut(X_4)$) is hyperbolic, but not convex-cocompact. 
\end{enumerate}
\end{eg}


\begin{eg}\label{eg-no(-2)curve}
Let $X$ be a K3 surface as in Theorem \ref{thm-geom-fin}, and suppose that it has no $(-2)$-curves. 
Then $\Aut(X)$ is a lattice by Theorem \ref{chara-lattice}. 
If $X$ is elliptic, 
then the Mordell--Weil rank of each elliptic fibration is $\rho_X-2$ by the theory of Kleinian groups (see \cite[Corollary 1 in Section 8 and Section 11.1]{Rat}). 
Especially when $\rho_X=4$, $\Aut(X)$ is hyperbolic and a lattice if and only if it is a uniform lattice. 



\end{eg}


\begin{eg}
Let $X$ be a K3 surface as in Theorem \ref{chara-conv-cocpt} (ii) with infinite automorphism group. 
By Theorem \ref{zero-entropy-classification}, 
when $\rho_X \geq 5$, $\Aut(X)$ is non-elementary relatively hyperbolic. 
%
For $\rho_X = 2, 3\text{ or }4$, 
there exists a K3 surface as in Theorem \ref{chara-conv-cocpt} (ii) whose automorphism group contains an infinite cyclic subgroup of finite index, hence elementary. 
This subgroup is generated by an automorphism of positive entropy (see \cite[Ch.15.2.5]{HuyK3}, \cite[Remark 3.10]{Brandhorst-Mezzedimi}, and Theorem \ref{zero-entropy-classification}).

\end{eg}


\begin{eg}\label{eg-vir-hyp-mfd-gp}
The following even hyperbolic lattices contain no roots and no nonzero isotropic vectors, and hence yield examples of complex K3 surfaces whose automorphism groups are virtually hyperbolic-manifold-groups as in Theorem \ref{K3-uniform-lattice} (ii): 
\[
\langle4\rangle\oplus\langle-8\rangle\oplus\langle-12k\rangle,
\qquad
\langle4\rangle\oplus\langle-8\rangle\oplus\langle-12\rangle\oplus\langle-12k\rangle, 
\]
where $k\in\ZZ_{>0}$ and $k\equiv1$ (mod $3$). 
\end{eg}





\begin{eg}
Let $X$ be a singular K3 surface over $\CC$. 
We denote by $\Hilb^n(X)$ the Hilbert scheme of $0$-dimensional closed subschemes of length $n$ of $X$. 
Then $\Aut(\Hilb^n(X))$ contains $F_2$ by Oguiso \cite[Theorem 1.3]{Oguiso-Tits}, 
hence is non-elementary relatively hyperbolic. 
\end{eg}


\bibliography{math}

@article {Amerik-Verbitsky,
    AUTHOR = {E. Amerik and M. Verbitsky},
     TITLE = {{Morrison--Kawamata cone conjecture for hyperk\"{a}hler manifolds}},
   JOURNAL = {Ann. Sci. \'{E}c. Norm. Sup\'{e}r.(4)},
  FJOURNAL = {},
    VOLUME = {50},
      YEAR = {2017},
    NUMBER = {4},
     PAGES = {973--993},
      ISSN = {},
   MRCLASS = {},
  MRNUMBER = {},
MRREVIEWER = {},
       DOI = {},
       URL = {},
}

@article {Bow1,
    AUTHOR = {B.H. Bowditch},
     TITLE = {{Geometrical Finiteness for Hyperbolic Groups}},
   JOURNAL = {J. Funct. Anal.},
  FJOURNAL = {},
    VOLUME = {113},
      YEAR = {1993},
    NUMBER = {},
     PAGES = {245--317},
      ISSN = {},
   MRCLASS = {},
  MRNUMBER = {},
MRREVIEWER = {},
       DOI = {},
       URL = {},
}

@article {Bow2,
    AUTHOR = {B.H. Bowditch},
     TITLE = {{Relatively hyperbolic groups}},
   JOURNAL = {Southampton Preprint},
  FJOURNAL = {},
    VOLUME = {},
      YEAR = {1999},
    NUMBER = {},
     PAGES = {},
      ISSN = {},
   MRCLASS = {},
  MRNUMBER = {},
MRREVIEWER = {},
       DOI = {},
       URL = {https://eprints.soton.ac.uk/29769/},
}

@article {Brandhorst-supersingular,
    AUTHOR = {S. Brandhorst},
     TITLE = {{Automorphisms of salem degree 22 on supersingular K3 surfaces of higher Artin invariant}},
   JOURNAL = {English. Math. Res. Lett.},
  FJOURNAL = {},
    VOLUME = {25},
      YEAR = {2018},
    NUMBER = {},
     PAGES = {1143--1150},
      ISSN = {},
   MRCLASS = {},
  MRNUMBER = {},
MRREVIEWER = {},
       DOI = {},
       URL = {},
}

@article {Brandhorst-Mezzedimi,
    AUTHOR = {S. Brandhorst and G. Mezzedimi},
     TITLE = {{Borcherds lattices and K3 surfaces of zero entropy}},
   JOURNAL = {to appear Amer. J. Math.},
  FJOURNAL = {},
    VOLUME = {preprint arXiv:2211.09600},
      YEAR = {2022},
    NUMBER = {},
     PAGES = {},
      ISSN = {},
   MRCLASS = {},
  MRNUMBER = {},
MRREVIEWER = {},
       DOI = {},
       URL = {},
}

@book {BriH,
    AUTHOR = {M.R. Bridson and A. Haefliger},
     TITLE = {{Metric Spaces of Non-Positive Curvature}},
    SERIES = {}, 
 PUBLISHER = {Springer-Verlag},
   ADDRESS = {},
      YEAR = {1999},
      ISBN = {},
}

@article {BLvL,
    AUTHOR = {M. Bright and A. Logan and R. v. Luijk},
     TITLE = {{Finiteness results for K3 surfaces over arbitrary fields}},
   JOURNAL = {Eur. J. Math.},
  FJOURNAL = {},
    VOLUME = {6},
      YEAR = {2020},
    NUMBER = {},
     PAGES = {336--366},
      ISSN = {},
   MRCLASS = {},
  MRNUMBER = {},
MRREVIEWER = {},
       DOI = {},
       URL = {},
}

@article {Cantat-K3,
    AUTHOR = {S. Cantat},
     TITLE = {{Dynamique des automorphismes des surfaces K3}},
   JOURNAL = {Acta Math.},
  FJOURNAL = {},
    VOLUME = {187},
      YEAR = {2001},
    NUMBER = {},
     PAGES = {1--57},
      ISSN = {},
   MRCLASS = {},
  MRNUMBER = {},
MRREVIEWER = {},
       DOI = {},
       URL = {},
}

@article {Cantat-Dolgachev,
    AUTHOR = {S. Cantat and I. V. Dolgachev},
     TITLE = {{Rational surfaces with a large group of automorphisms}},
   JOURNAL = {J. Amer. Math. Soc. },
  FJOURNAL = {},
    VOLUME = {25},
      YEAR = {2012},
    NUMBER = {},
     PAGES = {863--905},
      ISSN = {},
   MRCLASS = {},
  MRNUMBER = {},
MRREVIEWER = {},
       DOI = {},
       URL = {},
}

@article {DolgachevMartin2019,
    AUTHOR = {I. Dolgachev and G. Martin},
     TITLE = {{Numerically trivial automorphisms of {Enriques} surfaces in characteristic 2}},
   JOURNAL = {J. Math. Soc. Japan},
  FJOURNAL = {},
    VOLUME = {71},
      YEAR = {2019},
    NUMBER = {4},
     PAGES = {1181--1200},
      ISSN = {},
   MRCLASS = {},
  MRNUMBER = {},
MRREVIEWER = {},
       DOI = {},
       URL = {},
}

@article {Dolgachev-Martin,
    AUTHOR = {I. Dolgachev and G. Martin},
     TITLE = {{Automorphism groups of rational elliptic and quasi-elliptic surfaces in all characteristics}},
   JOURNAL = {Adv. Math.},
  FJOURNAL = {},
    VOLUME = {400},
      YEAR = {2022},
    NUMBER = {},
     PAGES = {108274},
      ISSN = {},
   MRCLASS = {},
  MRNUMBER = {},
MRREVIEWER = {},
       DOI = {},
       URL = {},
}

@article {Dolgachev-Zhang,
    AUTHOR = {I. V. Dolgachev and D.-Q. Zhang},
     TITLE = {{Coble rational surfaces}},
   JOURNAL = {Amer. J. Math.},
  FJOURNAL = {},
    VOLUME = {123},
      YEAR = {2001},
    NUMBER = {},
     PAGES = {79--114},
      ISSN = {},
   MRCLASS = {},
  MRNUMBER = {},
MRREVIEWER = {},
       DOI = {},
       URL = {},
}

@article {Gro,
    AUTHOR = {M. Gromov},
     TITLE = {{On the entropy of holomorphic maps}},
   JOURNAL = {Enseign. Math.},
  FJOURNAL = {},
    VOLUME = {49},
      YEAR = {2003},
    NUMBER = {},
     PAGES = {217--235},
      ISSN = {},
   MRCLASS = {},
  MRNUMBER = {},
MRREVIEWER = {},
       DOI = {},
       URL = {},
}

@book {HuyK3,
    AUTHOR = {D. Huybrechts},
     TITLE = {{Lectures on  K3 surfaces}},
    SERIES = {Cambridge Studies in Advanced Mathematics}, 
 PUBLISHER = {Cambridge University Press},
   ADDRESS = {},
      YEAR = {2016},
      ISBN = {},
}

@article {Kap,
    AUTHOR = {M. Kapovich},
     TITLE = {{Kleinian Groups in Higher Dimensions}},
   JOURNAL = {Geometry and Dynamics of Groups and Spaces},
  FJOURNAL = {},
    VOLUME = {},
      YEAR = {2007},
    NUMBER = {},
     PAGES = {487--564},
      ISSN = {},
   MRCLASS = {},
  MRNUMBER = {},
MRREVIEWER = {},
       DOI = {},
       URL = {},
}

@article {Kawamata,
    AUTHOR = {Y. Kawamata},
     TITLE = {{On the cone of divisors of Calabi–Yau fiber spaces}},
   JOURNAL = {Internat. J. Math.},
  FJOURNAL = {},
    VOLUME = {8},
      YEAR = {1997},
    NUMBER = {5},
     PAGES = {665--687},
      ISSN = {},
   MRCLASS = {},
  MRNUMBER = {},
MRREVIEWER = {},
       DOI = {},
       URL = {},
}

@article {KY,
    AUTHOR = {N. Kurnosov and E. Yasinsky},
     TITLE = {{Automorphisms of Hyperkähler Manifolds and Groups Acting on CAT(0) Spaces}},
   JOURNAL = {Birational Geometry, Kähler–Einstein Metrics and Degenerations},
  FJOURNAL = {},
    VOLUME = {},
      YEAR = {2023},
    NUMBER = {},
     PAGES = {477-500},
      ISSN = {},
   MRCLASS = {},
  MRNUMBER = {},
MRREVIEWER = {},
       DOI = {},
       URL = {},
}

@article {Lieblich-Maulik,
    AUTHOR = {M. Lieblich and D. Maulik},
     TITLE = {{A note on the cone conjecture for K3 surfaces in positive characteristic}},
   JOURNAL = {Math. Res. Lett.},
  FJOURNAL = {},
    VOLUME = {25},
      YEAR = {2018},
    NUMBER = {6},
     PAGES = {1879--1891},
      ISSN = {},
   MRCLASS = {},
  MRNUMBER = {},
MRREVIEWER = {},
       DOI = {},
       URL = {},
}

@article {Markman,
    AUTHOR = {E. Markman},
     TITLE = {{A survey of Torelli and monodromy results for holomorphic-symplectic varieties}},
   JOURNAL = {Complex and differential geometry, {\rm Springer}},
  FJOURNAL = {},
    VOLUME = {},
      YEAR = {2011},
    NUMBER = {},
     PAGES = {257--322},
      ISSN = {},
   MRCLASS = {},
  MRNUMBER = {},
MRREVIEWER = {},
       DOI = {},
       URL = {},
}

@article {MMV,
    AUTHOR = {G. Martin and G. Mezzedimi and D. C. Veniani
},
     TITLE = {{Enriques surfaces of zero entropy}},
   JOURNAL = {},
  FJOURNAL = {},
    VOLUME = {preprint arXiv:2406.18407},
      YEAR = {2024},
    NUMBER = {},
     PAGES = {},
      ISSN = {},
   MRCLASS = {},
  MRNUMBER = {},
MRREVIEWER = {},
       DOI = {},
       URL = {},
}

@article {Mezzedimi-19,
    AUTHOR = {G. Mezzedimi},
     TITLE = {{K3 surfaces of zero entropy admitting an elliptic fibration with only irreducible fibers}},
   JOURNAL = {J. Algebra},
  FJOURNAL = {},
    VOLUME = {587},
      YEAR = {2021},
    NUMBER = {},
     PAGES = {344--389},
      ISSN = {},
   MRCLASS = {},
  MRNUMBER = {},
MRREVIEWER = {},
       DOI = {},
       URL = {},
}

@article {MukConj,
    AUTHOR = {S. Mukai},
     TITLE = {{}},
   JOURNAL = {{\rm Proceedings of the 67th Algebra Symposium (Japanese)}},
  FJOURNAL = {},
    VOLUME = {},
      YEAR = {2018},
     PAGES = {},
      ISSN = {},
   MRCLASS = {},
  MRNUMBER = {},
MRREVIEWER = {},
       URL = {http://www.mathsoc.jp/assets/file/sections/algebra/algsympo/algsymp18/houkokusyu/11-Mukai.pdf},
}

@article {Namikawa,
    AUTHOR = {Y. Namikawa},
     TITLE = {{Periods of Enriques surfaces}},
   JOURNAL = {Math. Ann.},
  FJOURNAL = {},
    VOLUME = {270},
      YEAR = {1985},
    NUMBER = {},
     PAGES = {201--222},
      ISSN = {},
   MRCLASS = {},
  MRNUMBER = {},
MRREVIEWER = {},
       DOI = {},
       URL = {},
}

@article {Nikulin,
    AUTHOR = {V. V. Nikulin},
     TITLE = {{Finite automorphism groups of K\"{a}hler surfaces of type K3}},
   JOURNAL = {Proc. Moscow Math. Soc.},
  FJOURNAL = {},
    VOLUME = {38},
      YEAR = {1979},
    NUMBER = {},
     PAGES = {75--137},
      ISSN = {},
   MRCLASS = {},
  MRNUMBER = {},
MRREVIEWER = {},
       DOI = {},
       URL = {},
}

@article {Nikulin83,
    AUTHOR = {V. V. Nikulin},
     TITLE = {{Factor groups of groups of automorphisms of hyperbolic forms with respect to subgroups generated by 2-reflections. Algebrogeometric applications}},
   JOURNAL = {J. Math. Sci.},
  FJOURNAL = {},
    VOLUME = {22},
      YEAR = {1983},
    NUMBER = {4},
     PAGES = {1401--1475},
      ISSN = {},
   MRCLASS = {},
  MRNUMBER = {},
MRREVIEWER = {},
       DOI = {},
       URL = {},
}

@article {Oguiso-Tits,
    AUTHOR = {K. Oguiso},
     TITLE = {{Tits alternative in hypek\"{a}hler manifolds}},
   JOURNAL = {Math. Res. Lett.},
  FJOURNAL = {},
    VOLUME = {13},
      YEAR = {2006},
    NUMBER = {},
     PAGES = {307--316},
      ISSN = {},
   MRCLASS = {},
  MRNUMBER = {},
MRREVIEWER = {},
       DOI = {},
       URL = {},
}

@article {Oguiso-20-zero,
    AUTHOR = {K. Oguiso},
     TITLE = {{Automorphisms of hypek\"{a}hler manifolds in the view of topological entropy}},
   JOURNAL = {In: Contemp. Math.},
  FJOURNAL = {},
    VOLUME = {422},
      YEAR = {2007},
    NUMBER = {},
     PAGES = {173--185},
      ISSN = {},
   MRCLASS = {},
  MRNUMBER = {},
MRREVIEWER = {},
       DOI = {},
       URL = {},
}

@article {Oguiso-Sakurai,
    AUTHOR = {K. Oguiso and J. Sakurai},
     TITLE = {{Calabi–Yau threefolds of quotient type}},
   JOURNAL = {Asian J. Math.},
  FJOURNAL = {},
    VOLUME = {5},
      YEAR = {2001},
    NUMBER = {1},
     PAGES = {43--77},
      ISSN = {},
   MRCLASS = {},
  MRNUMBER = {},
MRREVIEWER = {},
       DOI = {},
       URL = {},
}

@article {Osin,
    AUTHOR = {D. Osin},
     TITLE = {{Relatively hyperbolic groups: intrinsic geometry, algebraic properties, and algorithmic problems}},
   JOURNAL = {Mem. Amer. Math. Soc.},
  FJOURNAL = {},
    VOLUME = {179},
      YEAR = {2006},
    NUMBER = {843},
     PAGES = {},
      ISSN = {},
   MRCLASS = {},
  MRNUMBER = {},
MRREVIEWER = {},
       DOI = {},
       URL = {},
}

@book {Rat,
    AUTHOR = {J.G. Ratcliffe},
     TITLE = {{Foundations of Hyperbolic Manifolds, 3rd ed.}},
    SERIES = {Graduate Texts in Mathematics}, 
 PUBLISHER = {Springer Nature},
   ADDRESS = {},
      YEAR = {2006},
      ISBN = {},
}

@article {Shi,
    AUTHOR = {T. Shioda},
     TITLE = {{Mordell--Weil Lattice of Higher Genus Fibration on a Fermat Surface}},
   JOURNAL = {J. Math. Sci. Univ. Tokyo},
  FJOURNAL = {},
    VOLUME = {22},
      YEAR = {2015},
    NUMBER = {},
     PAGES = {443--468},
      ISSN = {},
   MRCLASS = {},
  MRNUMBER = {},
MRREVIEWER = {},
       DOI = {},
       URL = {},
}

@book {SI,
    AUTHOR = {T. Shioda and H. Inose},
     TITLE = {{On singular K3 surfaces}},
    SERIES = {In Complex analysis and algebraic geometry, 119--136}, 
 PUBLISHER = {Iwanami Shoten},
   ADDRESS = {},
      YEAR = {1977},
      ISBN = {},
}

@article {Sterk,
    AUTHOR = {H. Sterk},
     TITLE = {{Finiteness results for algebraic K3 surfaces}},
   JOURNAL = {Math. Z.},
  FJOURNAL = {},
    VOLUME = {189},
      YEAR = {1985},
    NUMBER = {},
     PAGES = {507--514},
      ISSN = {},
   MRCLASS = {},
  MRNUMBER = {},
MRREVIEWER = {},
       DOI = {},
       URL = {},
}

@article {Takamatsu,
    AUTHOR = {T. Takamatsu},
     TITLE = {{On the finiteness of twists of irreducible symplectic varieties}},
   JOURNAL = {Math. Ann.},
  FJOURNAL = {},
    VOLUME = {392},
      YEAR = {2025},
    NUMBER = {1},
     PAGES = {339--371},
      ISSN = {},
   MRCLASS = {},
  MRNUMBER = {},
MRREVIEWER = {},
       DOI = {https://doi.org/10.1017/nmj.2025.10082},
       URL = {},
}

@article {Takatsu,
    AUTHOR = {T. Takatsu},
     TITLE = {{Blown-up boundaries associated with ample cones of K3 surfaces}},
   JOURNAL = {Nagoya Math. J.},
  FJOURNAL = {Nagoya Mathematical Journal},
    VOLUME = {260},
      YEAR = {2025},
    NUMBER = {},
     PAGES = {796--825},
      ISSN = {},
   MRCLASS = {},
  MRNUMBER = {},
MRREVIEWER = {},
       DOI = {https://doi.org/10.1017/nmj.2025.10082},
       URL = {},
}

@article {Tot,
    AUTHOR = {B. Totaro},
     TITLE = {{The cone conjecture for Calabi--Yau pairs in dimension two}},
   JOURNAL = {{\rm arXiv version of \cite{TotDuke}: arXiv:0901.3361}},
  FJOURNAL = {},
    VOLUME = {},
      YEAR = {2009},
    NUMBER = {},
     PAGES = {},
      ISSN = {},
   MRCLASS = {},
  MRNUMBER = {},
MRREVIEWER = {},
       DOI = {},
       URL = {},
}

@article {TotDuke,
    AUTHOR = {B. Totaro},
     TITLE = {{The cone conjecture for Calabi--Yau pairs in dimension 2}},
   JOURNAL = {Duke Math. J.},
  FJOURNAL = {},
    VOLUME = {154},
      YEAR = {2010},
    NUMBER = {2},
     PAGES = {241–263},
      ISSN = {},
   MRCLASS = {},
  MRNUMBER = {},
MRREVIEWER = {},
       DOI = {},
       URL = {},
}

@article {Totaro-hyp-geom,
    AUTHOR = {B. Totaro},
     TITLE = {{Algebraic surfaces and hyperbolic geometry, In: Current Developments in Algebraic Geometry}},
   JOURNAL = {Math. Sci. Res. Inst. Publ.},
  FJOURNAL = {},
    VOLUME = {59},
      YEAR = {2012},
    NUMBER = {},
     PAGES = {405--426},
      ISSN = {},
   MRCLASS = {},
  MRNUMBER = {},
MRREVIEWER = {},
       DOI = {},
       URL = {},
}

@article {Tuk,
    AUTHOR = {P. Tukia},
     TITLE = {{Convergence groups and Gromov’s metric hyperbolic spaces}},
   JOURNAL = {New Zealand J. Math.},
  FJOURNAL = {},
    VOLUME = {23},
      YEAR = {1994},
    NUMBER = {2},
     PAGES = {157--187},
      ISSN = {},
   MRCLASS = {},
  MRNUMBER = {},
MRREVIEWER = {},
       DOI = {},
       URL = {},
}

@article {Vin,
    AUTHOR = {\`{E}. Vinberg},
     TITLE = {{The two most algebraic K3 surfaces}},
   JOURNAL = {Math. Ann.},
  FJOURNAL = {},
    VOLUME = {265},
      YEAR = {1983},
    NUMBER = {1},
     PAGES = {1--21},
      ISSN = {},
   MRCLASS = {},
  MRNUMBER = {},
MRREVIEWER = {},
       DOI = {},
       URL = {},
}

@article {Wang,
    AUTHOR = {L. Wang},
     TITLE = {{On automorphisms and the cone conjecture for Enriques surfaces in odd characteristic}},
   JOURNAL = {Math. Res. Lett.},
  FJOURNAL = {},
    VOLUME = {28},
      YEAR = {2021},
    NUMBER = {4},
     PAGES = {1263--1281},
      ISSN = {},
   MRCLASS = {},
  MRNUMBER = {},
MRREVIEWER = {},
       DOI = {},
       URL = {},
}

@article {FulinXu,
    AUTHOR = {F. Xu},
     TITLE = {{On the cone conjecture for certain pairs of dimension at most $4$}},
   JOURNAL = {},
  FJOURNAL = {},
    VOLUME = {preprint arXiv:2405.20899},
      YEAR = {2024},
    NUMBER = {},
     PAGES = {},
      ISSN = {},
   MRCLASS = {},
  MRNUMBER = {},
MRREVIEWER = {},
       DOI = {},
       URL = {},
}

@article {Yom,
    AUTHOR = {Y. Yomdin},
     TITLE = {{Volume growth and entropy}},
   JOURNAL = {Israel J. Math. },
  FJOURNAL = {},
    VOLUME = {57},
      YEAR = {1987},
    NUMBER = {},
     PAGES = {285--300},
      ISSN = {},
   MRCLASS = {},
  MRNUMBER = {},
MRREVIEWER = {},
       DOI = {},
       URL = {},
}

@article {Yu-supersingular,
    AUTHOR = {X. Yu},
     TITLE = {{Elliptic fibrations on K3 surfaces and Salem numbers of maximal degree}},
   JOURNAL = {J. Math. Soc. Japan},
  FJOURNAL = {},
    VOLUME = {70},
      YEAR = {2018},
    NUMBER = {},
     PAGES = {1151--1163},
      ISSN = {},
   MRCLASS = {},
  MRNUMBER = {},
MRREVIEWER = {},
       DOI = {},
       URL = {},
}

@article {Yu,
    AUTHOR = {X. Yu},
     TITLE = {{K3 surface entropy and automorphism groups}},
   JOURNAL = {J. Algebraic Geom.},
  FJOURNAL = {},
    VOLUME = {34},
      YEAR = {2025},
    NUMBER = {},
     PAGES = {205--231},
      ISSN = {},
   MRCLASS = {},
  MRNUMBER = {},
MRREVIEWER = {},
       DOI = {},
       URL = {},
}
\bibliographystyle{alpha}
\end{document}